\documentclass[12pt,fleqn]{article}
\usepackage{amsmath,amsthm,amscd,amssymb}
\usepackage{latexsym}
\title{Concrete Constructions of Real Equiangular Line Sets}
\author{Janet C. Tremain\thanks{The author was supported by NSF DMS 0704216}\\Department of Mathematics\\University
of Missouri\\Columbia, MO 65211-4100\\\texttt{janet@math.missouri.edu}}

\newtheorem{Thm}{Theorem}[section]

\newtheorem{corollary}[Thm]{Corollary}

\newtheorem{proposition}[Thm]{Proposition}

\newtheorem{remark}[Thm]{Remark}

\newtheorem{theorem}[Thm]{Theorem}
\newtheorem{problem}[Thm]{Problem}

\def \R{\mathbb R}

\def \beq{\begin{eqnarray*}}
\def \eeq{\end{eqnarray*}}

\begin{document}

\footnotetext{\textit {Math Subject Classification.}  Primary: 51M04}

\maketitle

\begin{abstract}
We give some concrete constructions of real equiangular line sets.  The
emphasis here is on {\em building blocks} for certain angles which
are then used to build up larger equiangular line sets.  We concentrate on
angles greater than or equal to $1/7$.
\end{abstract}
\pagebreak

\tableofcontents

\pagebreak

\section{Index}
\setcounter{equation}{0}

\noindent {\bf In these notes,
we are constructing $M$ equiangular lines in
$\R^N$} at angle $1/\alpha$.
\vskip 12pt

\noindent (I)  N=2, M=3.
\vskip12pt

\noindent (II)  N=3, M=6.
\vskip12pt

\noindent (III)  {\bf Angle $1/3$}.
\vskip12pt

(A)  {\bf Building Blocks}
\vskip12pt

\ \ 1.  N=3, M=4.

\ \ 2.  N=4, M=6.

\ \ 3.  N=5, M=8.

\ \ 4.  N=7, M=8.

\vskip12pt
(B)  {\bf Specific Dimensions}
\vskip12pt

\ \ 1.  N=4, M=6.

\ \ 2.  N=5, M=8 

\ \ 3.  N=5, M=8 

\ \ 4.  N=5, M=10.

\ \ 5.  N=6, M=12.

\ \ 6.  N=6, M=12 (Another example).

\ \ 7.  N=6, M=16 (All entries $\pm 1/\sqrt{6}$).

\ \ 8.  N=6, M=16 (From building blocks).

\ \ 9.  N=7, M=28.

\  10.  N=7, M=28 (Simpler Format).

\ 11.  N=7, M=28 (Different Format).

\ 12.  N=7, M=28 (Another example).
\vskip12pt
(C)  {\bf General Sets}
\vskip12pt

\ \ 1.  $\R^N$ always has $M=2(N-1)$.  
\vskip12pt

\noindent (IV)  {\bf Angle $1/5$}.
\vskip12pt

(A)  {\bf Building Blocks}
\vskip12pt

\ \ 1.  N=4, M=4.

\ \ 2.  N=5, M=5 (Circulant containing $4 \times 4$ orthogonal).

\ \ 3.  N=5, M=5 (Circulant)

\ \ 4.  N=5, M=5 (Circulant)

\ \ 5.  N=7, M=4.

\ \ 6.  N=7, M=8.

\ \ 7.  N=8, M=8.

\vskip12pt

(B)  {\bf Specific Dimensions}
\vskip12pt

\ \ 1.  N=14, M=28.

\ \ 2.  N=15, M=30.

\ \ 3.  N=15, M=36.

\ \ 4.  N=16, M=40.

\ \ 5.  N=19, M= 60.

\ \ 6.  N=21, M=45.

\ \ 7.  N=21, M= 45 (Using building blocks)

\vskip12pt

(C)  {\bf General sets}
\vskip12pt

\ \ 1.  $\R^N$ always has $M=N-1$ vectors at angle $1/5$.

\ \ 2.  $\R^N$ always has $M=N-1$ vectors at angle $1/5$ (Another example).

\ \ 3.  $\R^{3N+1}$ always has $4N$ vectors at angle $1/5$.
\vskip12pt

\noindent (V)  {\bf Angle $1/7$}.
\vskip12pt

(A)  {\bf Building Blocks}
\vskip12pt

\ \ 1.  N=7, M=7.

\ \ 2.  N=7, M=8.

\vskip12pt

\noindent (IX)  {\bf $M=2N$ vectors in $\R^N$}.
\vskip12pt

1.  N=3, M=6, angle $1/\sqrt{5}$.

2.  N=5, M=10, angle $1/3$.

3.  N=6, M=12, angle $1/3$.

3.  N=14, M=28, angle $1/5$.

4.  N=15, M=30, angle $1/5$.

\vskip12pt

\noindent (X)  {\bf $M=N+1$ vectors in $\R^N$}
\vskip12pt

\noindent (XI)  {\bf Unitary Matrices}
\vskip12pt

1.  $M=N$ (Circulant self-adjoint unitary with two entries)

2.  $M=N$ (Circulant self-adjoint unitary with two entries)

3.  $M=N$ (Circulant unitary with two entries)
\vskip12pt

\noindent (XII)  {\bf Multiple Angles}
\vskip12pt

5.  $M=2N$, inner products $\pm 1/\sqrt{5},0$.
(Made up of two circulants).

\newpage

\section{Introduction}

A very old problem is:  

\begin{problem}  How many equiangular lines can we draw
through the origin in $\R^N$?
\end{problem}

  This means, if we choose a set of
unit length vectors $\{f_m\}_{m=1}^M$, one on each line, then 
there is a constant $c$ so that for
all $1\le m\not= n \le M$ we have
\[ |\langle f_n,f_m\rangle | = c.\]
These inner products
represent the cosine of the acute angle between the lines.
The problem of constructing any number (especially, the maximal number)
of equiangular lines in ${\mathbf R}^N$ is one of the most elementary
and at the same time one of the most difficult problems in mathematics.
After sixty years of research, the maximal number of equiangular lines
in ${\mathbf R}^N$ is known only for 35 dimensions.  
For a slightly more general view of this topic see Benedetto and
Kolesar \cite{BK}.  This line of
research was started in 1948 by Hanntjes \cite{H} in the setting of elliptic
geometry where he identified the maximal
number of equiangular lines in ${\mathbf R}^N$
for $n=2,3$.  Later, Van Lint and Seidel \cite{LiS} classified the
largest number of equiangular lines in ${\mathbf R}^N$ for dimensions
$N\le 7$ and at the same time emphasized the relations to
discrete mathematics.  
In 1973, Lemmens and Seidel \cite{LS} made a comprehensive
study of real equiangular line sets which is still today a fundamental
piece of work.  Gerzon \cite{LS} gave an upper bound for the maximal
number of equiangular lines in ${\mathbf R}^N$:

\begin{theorem}[Gerzon] \label{T2}If we have $M$ 
equiangular lines in ${\mathbf R}^N$ then
$$M\le \frac{N(N+1)}{2}.$$
\end{theorem}
\vskip8pt
We will see that in most cases there are many fewer lines than this bound gives.
  Also, P. Neumann \cite{LS} produced a fundamental result
in the area:

\begin{theorem}[P. Neumann]\label{T1}
If ${\mathbf R}^N$ has
$M$ equiangular lines at angle $1/\alpha$ and $M>2N$, then $\alpha$ is an odd
integer.  
\end{theorem}
\vskip8pt
Finally, there is a lower bound on the angle formed by equiangular
line sets.  
\begin{theorem}\label{T5}
If $\{f_m\}_{m=1}^M$ is a set of norm one vectors in 
${\mathbf R}^N$, then
\begin{equation}\label{E1}
\max_{m\not= n}|\langle f_m,f_n\rangle| \ge \sqrt{\frac{M-N}{N(M-1)}}.
\end{equation}
Moreover, we have equality if and only if $\{f_m\}_{m=1}^M$ is an
equiangular tight frame and in this case the tight frame bound
is $\frac{M}{N}$.      
\end{theorem}
\vskip8pt

This inequality goes back to Welch \cite{W}.  Strohmer
and Heath \cite{SH} and Holmes and Paulsen
\cite{HP} give more direct arguments which also yields
the "moreover" part.  For some reason, in the literature there is a further
assumption added to the "moreover" part of Theorem \ref{T5} that
the vectors span ${\mathbf R}^N$.  This assumption is not necessary.
That is, equality in inequality \ref{E1} already implies that the
vectors span the space \cite{CRT}.  

The status of the equiangular line problem at this point is
summarized in the following chart \cite{LS,CRT,CRT1} where $N$ is the dimension
of the Hilbert space, $M$ is the maximal number of equiangular lines
and these will occur at the angle $1/\alpha$.
\vskip12pt

\noindent {\bf Table I:  Maximal equiangular line sets}

\begin{equation*}\begin{split}&\begin{array}{|rccccccccc|}\hline
N= & 2 & 3 & 4 & 5 & 6 & 7 & \dots & 13 & 14\\
M= & 3 & 6 & 6 & 10 & 16 & 28 & \dots & 28&28-30\\
\alpha= & 2 & \sqrt{5} & 3 & 3 & 3 & 3 & \dots &3&5\\ \hline\end{array}\\[6pt]
&\begin{array}{|rcccccc|}\hline
N= & 15 & 16 & 17 & 18 & 19 & 20 \\
M= & 36 & \ge 40 & \ge 48 & \ge 48 & 72-76* & 92-96* \\
\alpha= & 5 & 5 & 5 & 5 & 5 & 5 \\ \hline\end{array}\\[6pt]
&\begin{array}{|rccccccc|}\hline
N= & 21 & 22& 23 & \dots & 41 & 42&43\\
M= & 126 & 176& 276 & \dots & 276 & \ge 276 &344\\
\alpha= & 5 & 5& 5 & \dots & 5 & 5&7\\ \hline\end{array}\end{split}\end{equation*}

The $*$ in the chart represents two cases which have been
reported as solved in the literature but actually are still open.

For recent results on the equiangular line problem see \cite{CRT,CRT1}.

Here we will give concrete constructions for real equiangular line
sets.  The emphasis will be on constructing {\bf building blocks}
which can be put together to build larger and larger equiangular
line sets.

\newpage

\section{$M=3$ vectors in $\R^2$ at angle $\frac{1}{2}$}

\[ \begin{bmatrix} f_1 \\ f_2 \\ f_3
\end{bmatrix} =
\begin{bmatrix} 0 & 1 \\
-\frac{\sqrt{3}}{2} & - \frac{1}{2}\\
\frac{\sqrt{3}}{2} & - \frac{1}{2}
\end{bmatrix}
\]

\section{$M=6$ vectors in $\R^3$ at angle $\frac{1}{\sqrt{5}}$}
\[
\begin{bmatrix}
0 & \sqrt{\frac{5-\sqrt{5}}{10}} & \sqrt{\frac{5+\sqrt{5}}{10}}\\
0 & -\sqrt{\frac{5-\sqrt{5}}{10}} & \sqrt{\frac{5+\sqrt{5}}{10}}\\
\sqrt{\frac{5-\sqrt{5}}{10}} & \sqrt{\frac{5+\sqrt{5}}{10}} & 0\\
-\sqrt{\frac{5-\sqrt{5}}{10}} & \sqrt{\frac{5+\sqrt{5}}{10}} & 0\\
\sqrt{\frac{5+\sqrt{5}}{10}} & 0 & \sqrt{\frac{5-\sqrt{5}}{10}}\\
\sqrt{\frac{5+\sqrt{5}}{10}} & 0 & -\sqrt{\frac{5-\sqrt{5}}{10}}
\end{bmatrix}\]

\pagebreak

\section{Angle $1/3$}

\noindent {\bf (A)  Building Blocks}
\vskip12pt

\subsection{$M=4$ vectors in $\R^3$}
\[ \frac{1}{\sqrt{3}}\begin{bmatrix}
1&1&1\\
-1&1&1\\
1&-1&1\\
1&1&-1
\end{bmatrix}\]

\subsection{$M=6$ vectors in $\R^4$ at angle $1/3$}

\[ \begin{bmatrix}
\sqrt{\frac{1}{3}}& \sqrt{\frac{2}{3}}&&\\
\sqrt{\frac{1}{3}}& -\sqrt{\frac{2}{3}}&&\\
\sqrt{\frac{1}{3}}&& \sqrt{\frac{2}{3}}&\\
\sqrt{\frac{1}{3}}&& -\sqrt{\frac{2}{3}}&\\
\sqrt{\frac{1}{3}}&&& \sqrt{\frac{2}{3}}\\
\sqrt{\frac{1}{3}}&&& -\sqrt{\frac{2}{3}}\\
\end{bmatrix}\]

\subsection{$M=8$ vectors in $\R^5$ at angle $1/3$}

\[ \sqrt{\frac{1}{3}}\begin{bmatrix}
1 & 1 & 1 & 0 & 0 \\
-1 & 1 & 1 & 0 & 0\\
1 & -1 & 1 & 0 & 0\\
1 & 1 & -1 & 0 & 0 \\
0 & 0 & 1 & 1 & 1\\
0 & 0 & -1 & 1 & 1\\
0 & 0 & 1 & -1 & 1\\
0 & 0 & 1 & 1 & -1
\end{bmatrix}\]

\subsection{$M=8$ vectors in $\R^7$ at angle $1/3$}

\[ \sqrt{\frac{1}{10}}\begin{bmatrix}
1&1&1&1&\sqrt{2}&\sqrt{2}& \sqrt{2}\\
1&1&-1&-1&\sqrt{2}&\sqrt{2}& -\sqrt{2}\\
1&-1&1&-1&\sqrt{2}&-\sqrt{2}& -\sqrt{2}\\
1&-1&-1&1&-\sqrt{2}&\sqrt{2}& -\sqrt{2}\\
1&1&1&1&-\sqrt{2}&-\sqrt{2}& -\sqrt{2}\\
1&1&-1&-1&-\sqrt{2}&-\sqrt{2}& \sqrt{2}\\
1&-1&1&-1&-\sqrt{2}&\sqrt{2}& \sqrt{2}\\
1&-1&-1&1&\sqrt{2}&-\sqrt{2}& \sqrt{2}\\
\end{bmatrix}\]

\pagebreak

\noindent {\bf (B)  Specific Dimensions}
\vskip12pt

\subsection{$M=6$ vectors in $\R^4$ at angle $1/3$}

\[
\begin{bmatrix}
1 & 0 & 0 & 0\\
-\frac{1}{3} & \frac{1}{3}\sqrt{\frac{3}{2}} & \frac{1}{3\sqrt{2}}
& \sqrt{\frac{2}{3}}\\
- \frac{1}{3} & \frac{1}{3}\sqrt{\frac{3}{2}} & \frac{1}{3\sqrt{2}}
& -\sqrt{\frac{2}{3}}\\
-\frac{1}{3} & 0 & \frac{2\sqrt{2}}{3} & 0\\
- \frac{1}{3} & - \sqrt{\frac{2}{3}} & - \frac{\sqrt{2}}{3} & 0 \\
- \frac{1}{3} & \sqrt{\frac{2}{3}} &- \frac{\sqrt{2}}{3} & 0
\end{bmatrix}\]

\begin{remark}
The above example shows that
 an $M$ element equiangular tight frame for $\R^N$ need not have the property
that every $N$ element subset is linearly independent?  In the above
example, we clearly have four vectors which sit in $\R^3$.
\end{remark}
\pagebreak

\subsection{$M=8$ vectors in $\R^5$ at angle $1/3$}

\begin{equation*}
\begin{aligned}[b] &\quad 
&\left[ \begin{array}{ccccc}
\frac{1}{\sqrt{6}}& \frac{1}{\sqrt{6}}&\frac{1}{\sqrt{6}}
&\frac{1}{\sqrt{6}}&\frac{1}{\sqrt{3}}\\
\frac{1}{\sqrt{6}}& \frac{1}{\sqrt{6}}&\frac{-1}{\sqrt{6}}
&\frac{-1}{\sqrt{6}}&\frac{1}{\sqrt{3}}\\
\frac{1}{\sqrt{6}}& \frac{-1}{\sqrt{6}}&\frac{1}{\sqrt{6}}
&\frac{-1}{\sqrt{6}}&\frac{1}{\sqrt{3}}\\
\frac{1}{\sqrt{6}}& \frac{-1}{\sqrt{6}}&\frac{-1}{\sqrt{6}}
&\frac{1}{\sqrt{6}}&\frac{1}{\sqrt{3}}\\
\frac{1}{\sqrt{6}}& \frac{1}{\sqrt{6}}&\frac{1}{\sqrt{6}}
&\frac{1}{\sqrt{6}}&\frac{-1}{\sqrt{3}}\\
\frac{1}{\sqrt{6}}& \frac{1}{\sqrt{6}}&\frac{-1}{\sqrt{6}}
&\frac{-1}{\sqrt{6}}&\frac{-1}{\sqrt{3}}\\
\frac{1}{\sqrt{6}}& \frac{-1}{\sqrt{6}}&\frac{1}{\sqrt{6}}
&\frac{-1}{\sqrt{6}}&\frac{-1}{\sqrt{3}}\\
\frac{1}{\sqrt{6}}& \frac{-1}{\sqrt{6}}&\frac{-1}{\sqrt{6}}
&\frac{1}{\sqrt{6}}&\frac{-1}{\sqrt{3}}\\
\end{array}\right] \end{aligned}\]

\subsection{$M=8$ vectors in $\R^5$ at angle $1/3$}

\[ \sqrt{\frac{1}{3}}\begin{bmatrix}
1&1&1&0&0\\
-1&1&1&0&0\\
1&-1&1&0&0\\
1&1&-1&0&0\\
0&0&1&1&1\\
0&0&-1&1&1\\
0&0&1&-1&1\\
0&0&1&1&-1\\
\end{bmatrix}\]

\subsection{$M=10$ vectors in $\R^5$ at angle $1/3$}

\[\begin{bmatrix}
\sqrt{\frac{1}{3}} & \sqrt{\frac{2}{3}} & 0 & 0 & 0\\
-\sqrt{\frac{1}{3}} & \sqrt{\frac{2}{3}} & 0 & 0 & 0\\
\sqrt{\frac{1}{3}} & 0 & \sqrt{\frac{2}{3}} & 0 & 0\\
- \sqrt{\frac{1}{3}} & 0 &\sqrt{\frac{2}{3}} & 0 & 0\\
\sqrt{\frac{1}{3}} & 0 & 0 & \sqrt{\frac{2}{3}} & 0\\
-\sqrt{\frac{1}{3}} & 0 & 0 &\sqrt{\frac{2}{3}} & 0\\
0 & \frac{1}{3}\sqrt{\frac{3}{2}} & \frac{1}{3} \sqrt{\frac{3}{2}} &
\frac{1}{3}\sqrt{\frac{3}{2}} & \sqrt{\frac{1}{2}}\\
0 & \frac{1}{3}\sqrt{\frac{3}{2}} & -\frac{1}{3} \sqrt{\frac{3}{2}} &
\frac{1}{3}\sqrt{\frac{3}{2}} & -\sqrt{\frac{1}{2}}\\
0 & -\frac{1}{3}\sqrt{\frac{3}{2}} & -\frac{1}{3} \sqrt{\frac{3}{2}} &
\frac{1}{3}\sqrt{\frac{3}{2}} & \sqrt{\frac{1}{2}}\\
0 & \frac{1}{3}\sqrt{\frac{3}{2}} & -\frac{1}{3} \sqrt{\frac{3}{2}} &
-\frac{1}{3}\sqrt{\frac{3}{2}} & \sqrt{\frac{1}{2}}\\
\end{bmatrix}\]

\pagebreak

\subsection{$M=12$ vectors in $\R^6$ at angle $1/3$}

\begin{equation*}
\begin{aligned}[b] &\quad 
\left[ \begin{array}{cccccc}
1&2&3&4&5&6\\
\sqrt{\frac{1}{3}}& \sqrt{\frac{2}{3}}&&&&\\
\sqrt{\frac{1}{3}}&-\sqrt{\frac{2}{3}}&&&&\\
\sqrt{\frac{1}{3}} & & \sqrt{\frac{2}{3}}&&&\\
\sqrt{\frac{1}{3}}& & -\sqrt{\frac{2}{3}}&&&\\
\sqrt{\frac{1}{3}}& & & \sqrt{\frac{2}{3}}&&\\
\sqrt{\frac{1}{3}}& & & -\sqrt{\frac{2}{3}}&&\\
\sqrt{\frac{1}{3}}& & & & \sqrt{\frac{2}{3}}&\\
\sqrt{\frac{1}{3}}& & & &-\sqrt{\frac{2}{3}}&\\
& \frac{1}{\sqrt{6}}& \frac{1}{\sqrt{6}}&\frac{1}{\sqrt{6}}
&\frac{1}{\sqrt{6}}&\frac{1}{\sqrt{3}}\\
&\frac{1}{\sqrt{6}}& \frac{1}{\sqrt{6}}&\frac{-1}{\sqrt{6}}
&\frac{-1}{\sqrt{6}}&\frac{1}{\sqrt{3}}\\
&\frac{1}{\sqrt{6}}& \frac{-1}{\sqrt{6}}&\frac{1}{\sqrt{6}}
&\frac{-1}{\sqrt{6}}&\frac{1}{\sqrt{3}}\\
&\frac{1}{\sqrt{6}}& \frac{-1}{\sqrt{6}}&\frac{-1}{\sqrt{6}}
&\frac{1}{\sqrt{6}}&\frac{1}{\sqrt{3}}\\
\end{array}\right] \end{aligned}\]

\subsection{$M=12$ vectors in $\R^6$ at angle $1/3$ - another example}

We use a building block:
\[ \frac{1}{\sqrt{3}}\begin{bmatrix}
1&1&1\\
-1&1&1\\
1&-1&1\\
1&1&-1
\end{bmatrix}\]

Using this we have (using + for $1/\sqrt{3}$ and - for $-1/\sqrt{3}$):

\[ \begin{bmatrix}
& 1&2&3&4&5&6\\
1&+&+&+&&&\\
2&-&+&+&&&&\\
3&+&-&+&&&&\\
4&+&+&-&&&&\\
5&&&+&+&+&\\
6&&&-&+&+&\\
7&&&+&-&+&\\
8&&&+&+&-&\\
9&+&&&&+&+\\
10&-&&&&+&+\\
11&+&&&&-&+\\
12&+&&&&+&-
\end{bmatrix} \]

\pagebreak

\subsection{$M=16$ vectors in $\R^6$ at angle $1/3$}

\[ \frac{1}{\sqrt{6}}\begin{bmatrix}
1 & 1 & 1 & 1 & 1 & 1 \\
-1 & -1 & 1 & 1 & 1 & 1 \\
 -1 & 1 & -1 & 1 & 1 & 1\\
 -1 & 1 & 1 & -1 & 1 & 1  \\
 -1 & 1 & 1 & 1 & -1 & 1 \\
 -1 & 1 & 1 & 1 & 1 & -1\\
 1 & -1 & -1 & 1 & 1 & 1  \\
 1 & -1 & 1 & -1 & 1 & 1 \\
 1 & -1 & 1 & 1 & -1 & 1 \\
 1 & -1 & 1 & 1 & 1 & -1\\
 1 & 1 & -1 & -1 & 1 & 1 \\
 1 & 1 & -1 & 1 & -1 & 1\\
 1 & 1 & -1 & 1 & 1 & -1\\
 1 & 1 & 1 & -1 & -1 & 1\\
 1 & 1 & 1 & -1 & 1 & -1\\
 1 & 1 & 1 & 1 & -1 & -1
\end{bmatrix}
\]

\begin{proof}
The above matrix has mutual inner products $\pm 1/3$ since
each row $2\le m\le 16$ has two $-1's$.  So their inner
product with row 1 are all $1/3$.  The other 15 rows
come from putting $-1's$ in the 15 positions coming
from 16 choose 2.  So any two rows either have no
$-1's$ in common and their inner product is $-1/3$ or
they have one $-1$ in common and their inner product is
$1/3$.
\end{proof}

\begin{remark}This family is part of the $M=28$ equiangular
lines in $\R^7$ at angle $1/3$.
\end{remark}
\pagebreak

\subsection{$M=16$ vectors in $\R^6$ at angle $1/3$ from
building blocks}

\[ \sqrt{\frac{1}{3}}\begin{bmatrix}
1 & 1 & 1 & 0 & 0 & 0 \\
-1 & 1 & 1 & 0 & 0 & 0\\
1 & -1 & 1 & 0 & 0 & 0\\
1 & 1 & -1 & 0 & 0 & 0\\
0 & 0 & 1 & 1 & 1 & 0 \\
0 & 0 & -1 & 1 & 1 & 0\\
0 & 0 & 1 & -1 & 1 & 0\\
0 & 0 & 1 & 1 & -1 & 0\\
1 & 0 & 0 & 0 & 1 & 1\\
-1 & 0 & 0 & 0 & 1 & 1\\
1 & 0 & 0 & 0 & -1 & 1\\
1 & 0 & 0 & 0 & 1 & -1\\
0 & 1 & 0 & 1 & 0 & 1\\
0 & -1 & 0 & 1 & 0 & 1\\
0 & 1 & 0 & 1 & 0 & 1\\
0 & 1 & 0 & 1 & 0 & -1
\end{bmatrix}\]

\pagebreak

\subsection{$M=28$ vectors in $\R^7$ at angle
$1/3$}.

\[ \sqrt{\frac{1}{3}}\begin{bmatrix}
1 & 1 & 1 & 0 & 0 & 0 & 0 \\
-1 & 1 & 1 & 0 & 0 & 0 & 0\\
1 & -1 & 1 & 0 & 0 & 0 & 0\\
1 & 1 & -1 & 0 & 0 & 0 & 0\\
0 & 0 & 1 & 1 & 1 & 0 & 0 \\
0 & 0 & -1 & 1 & 1 & 0 & 0  \\
0 & 0 & 1 & -1 & 1 & 0 & 0\\
0 & 0 & 1 & 1 & -1 & 0 & 0\\
1 & 0 & 0 & 0 & 1 & 1 & 0\\
-1 & 0 & 0 & 0 & 1 & 1 & 0\\
1 & 0 & 0 & 0 &  -1 & 1 & 0\\
1 & 0 & 0 & 0 & 1 & -1 & 0\\
0 & 1 & 0 & 1 & 0 & 1 & 0\\
0 & -1 & 0 & 1 & 0 & 1 & 0\\
0 & 1 & 0 & -1 & 0 & 1 & 0\\
0 & 1 & 0 & 1 & 0 & -1 & 0\\
0 & 0 & 1 & 0 & 0 & 1 & 1 \\
0 & 0 & -1 & 0 & 0 & 1 & 1\\
0 & 0 & 1 & 0 & 0 & -1 & 1\\
0 & 0 & 1 & 0 & 0 & 1 & -1\\
1 & 0 & 0 & 1 & 0 & 0 & 1\\
-1 & 0 & 0 & 1 & 0 & 0 & 1\\
1 & 0 & 0 & -1 & 0 & 0 & 1\\
1 & 0 & 0 & 1 & 0 & 0 & -1\\
0 & 1 & 0 & 0 & 1 & 0 & 1\\
0 & -1 & 0 & 0 & 1 & 0 & 1 \\
0 & 1 & 0 & 0 & -1 & 0 & 1\\
0 & 1 & 0 & 0 & 1 & 0 & -1
\end{bmatrix}\]

\pagebreak

\subsection{$M=28$ vectors in $\R^7$ at angle $1/3$ - simpler
representation}
Let us try this without the zeroes to see if
it is clearer.  Also, we add the dimension numbers.

\begin{equation*}\begin{aligned} \sqrt{\frac{1}{3}}\end{aligned}
\begin{aligned}&\quad \begin{array}{p{.5cm}p{.5cm}p{.5cm}p{.5cm}p{.5cm}p{.5cm}p{.5cm}}
1 & 2 & 3 & 4 & 5 & 6 & 7\end{array}\\
&\left[ \begin{array}{p{.5cm}p{.5cm}p{.5cm}p{.5cm}p{.5cm}p{.5cm}p{.5cm}}
1 & 1 & 1 & & & &\\
-1 & 1 & 1 & & & &\\
1 & -1 & 1 & & & &\\
1 & 1 & -1 & & & &\\
& & 1 & 1 & 1 & &\\
& & -1 & 1 & 1 & &\\
& & 1 & -1 & 1 & &\\
& & 1 & 1 & -1 & &\\
1 & & & & 1 & 1 &\\
-1 & & & & 1 & 1 &\\
1 & & & & -1 & 1 & \\
1 & & & & 1 & -1 &\\
& 1 & & 1 & & 1 & \\
& -1 && 1 && 1 &\\
& 1 && -1 && 1 &\\
& 1 && 1 && -1 &\\
&& 1 &&& 1 & 1\\
&& -1 &&& 1 & 1\\
&& 1 &&& -1 & 1\\
&& 1 &&& 1 & -1\\
1 & & & 1 & & & 1\\
-1 & & & 1 & & & 1\\
1 & & & -1 & & & 1\\
1 & & & 1 & & & -1\\
& 1 & & & 1 & & 1\\
& -1 & & & 1 & & 1\\
& 1 & & & -1 & & 1\\
& 1 & & & 1 & & -1\\
\end{array}\right]\end{aligned}\end{equation*}

\pagebreak

\subsection{$M=28$ vectors in $\R^7$ at angle $1/3$ - Different format}

Yet another way to simplify this:

First, we use our {\em building block} for $4$ vectors
in $\R^3$ at angle $1/3$ where $+$ means $1$ and $-$ means
$-1$.

\[ BB1 =\sqrt{\frac{1}{3}}\begin{bmatrix}
+ & + & +\\
- & + & + \\
+ & - & + \\
+ & + & - \\
\end{bmatrix}\]
Now, we put these groups of four vectors into
our chart by just putting a {\em dot} where the
columns go.

\begin{equation*}\begin{aligned}[b]\begin{array}{c}
BB1\  4 \\
BB1\  4\\
BB1\  4 \\
BB1\ 4\\
BB1\ 4 \\
BB1\ 4\\
BB1\ 4 \\
\end{array}\end{aligned}
\begin{aligned}[b] &\quad \begin{array}{ccccccc}
1 & 2 & 3 & 4 & 5 & 6 & 7\end{array}\\
& \left[ \begin{array}{ccccccc}
\bullet & \bullet & \bullet & &&&\\
& & \bullet & \bullet & \bullet & &\\
\bullet & & & &\bullet & \bullet &\\
& \bullet && \bullet && \bullet &\\
& & \bullet & & & \bullet & \bullet\\
\bullet & & & \bullet & & &\bullet\\
& \bullet & & & \bullet & & \bullet\\
\end{array}\right]\end{aligned}
\end{equation*}

\pagebreak

\subsection{$M=28$ vectors in $\R^{7}$ at angle $1/3$ - second example}

We make a bulding block by embedding our $12$ vectors in $\R^6$
into $\R^7$.

\[ BB9=\frac{1}{\sqrt{6}}\begin{bmatrix}
0&1 & 1 & 1 & 1 & 1 & 1 \\
0&-1 & -1 & 1 & 1 & 1 & 1 \\
0& -1 & 1 & -1 & 1 & 1 & 1\\
0& -1 & 1 & 1 & -1 & 1 & 1  \\
0& -1 & 1 & 1 & 1 & -1 & 1 \\
0& -1 & 1 & 1 & 1 & 1 & -1\\
0& 1 & -1 & -1 & 1 & 1 & 1  \\
0& 1 & -1 & 1 & -1 & 1 & 1 \\
0& 1 & -1 & 1 & 1 & -1 & 1 \\
0& 1 & -1 & 1 & 1 & 1 & -1\\
0& 1 & 1 & -1 & -1 & 1 & 1 \\
0& 1 & 1 & -1 & 1 & -1 & 1\\
0& 1 & 1 & -1 & 1 & 1 & -1\\
0& 1 & 1 & 1 & -1 & -1 & 1\\
0& 1 & 1 & 1 & -1 & 1 & -1\\
0& 1 & 1 & 1 & 1 & -1 & -1
\end{bmatrix}
\]

Under this building block we put our standard $12$ vectors in $\R^7$
at angle $1/3$

\[ \begin{bmatrix}
1&2&3&4&5&6&7\\
\sqrt{\frac{1}{3}}&\sqrt{\frac{2}{3}}&&&&&\\
\sqrt{\frac{1}{3}}&-\sqrt{\frac{2}{3}}&&&&&\\
\sqrt{\frac{1}{3}}&&\sqrt{\frac{2}{3}}&&&&\\
\sqrt{\frac{1}{3}}&&-\sqrt{\frac{2}{3}}&&&&\\
\sqrt{\frac{1}{3}}&&&\sqrt{\frac{2}{3}}&&&\\
\sqrt{\frac{1}{3}}&&&-\sqrt{\frac{2}{3}}&&&\\
\sqrt{\frac{1}{3}}&&&&\sqrt{\frac{2}{3}}&&\\
\sqrt{\frac{1}{3}}&&&&-\sqrt{\frac{2}{3}}&&\\
\sqrt{\frac{1}{3}}&&&&&\sqrt{\frac{2}{3}}&\\
\sqrt{\frac{1}{3}}&&&&&-\sqrt{\frac{2}{3}}&\\
\sqrt{\frac{1}{3}}&&&&&&\sqrt{\frac{2}{3}}\\
\sqrt{\frac{1}{3}}&&&&&&-\sqrt{\frac{2}{3}}\\
\end{bmatrix} \]

\pagebreak
\noindent {\bf (C)  General Sets}
\vskip12pt

\section{$M=2(N-1)$ vectors in $\R^N$ at angle $1/3$}

\begin{proposition}\label{propos1}
For every $N$, $\R^N$ has $2(N-1)$
equiangular lines spanning $\R^N$ at angle
$1/3$.  This family is never a tight frame.  
\end{proposition}

\begin{proof}
We will just write down the vectors:
\[
\begin{bmatrix}
\sqrt{\frac{1}{3}} & \sqrt{\frac{2}{3}} & 0 & 0 & \cdots &
0 & 0 \\
\sqrt{\frac{1}{3}} & - \sqrt{\frac{2}{3}} & 0 & 0 & \cdots &
0 & 0 \\
\sqrt{\frac{1}{3}} & 0 & \sqrt{\frac{2}{3}} & 0 & \cdots & 0 & 0\\
\sqrt{\frac{1}{3}} & 0 & -\sqrt{\frac{2}{3}} & 0 & \cdots & 0 & 0\\
\vdots & \vdots & \vdots & \vdots & \cdots & 0 & 0\\
\sqrt{\frac{1}{3}} & 0 & 0 & 0 & \cdots & 0 & \sqrt{\frac{2}{3}}\\
\sqrt{\frac{1}{3}} & 0 & 0 & 0 & \cdots & 0 & -\sqrt{\frac{2}{3}}.
\end{bmatrix}\]
\end{proof}

\pagebreak

\section{Angle $1/5$}
\noindent {\bf (A)  Building Blocks}
\vskip12pt

\subsection{$M=4$ vectors in $\R^4$ at angle $1/5$}

\[  \begin{bmatrix}
\sqrt{\frac{2}{5}}& \sqrt{\frac{1}{5}}& \sqrt{\frac{1}{5}}
& \sqrt{\frac{1}{5}}\\
\sqrt{\frac{2}{5}}& -\sqrt{\frac{1}{5}}& -\sqrt{\frac{1}{5}}
& \sqrt{\frac{1}{5}}\\
\sqrt{\frac{2}{5}}& -\sqrt{\frac{1}{5}}& \sqrt{\frac{1}{5}}
& -\sqrt{\frac{1}{5}}\\
\sqrt{\frac{2}{5}}& \sqrt{\frac{1}{5}}& -\sqrt{\frac{1}{5}}
& -\sqrt{\frac{1}{5}}
\end{bmatrix}\]

\subsection{$M=5$ vectors in $\R^5$ at angle $1/5$}

\[ \sqrt{\frac{1}{5}}\begin{bmatrix}
1&1&1&1&1\\
-1&1&1&-1&1\\
-1&-1&1&1&-1\\
1&-1&1&-1&1\\
-1&-1&-1&1&1\\
\end{bmatrix}\]

\begin{remark}
Note that the $4 \times 4$ submatrix in the upper left corner
is an orthogonal matrix.
\end{remark}

\subsection{$M=5$ vectors in $\R^5$ at angle $1/5$
and is circulant}

\[ 
 \frac{1}{\sqrt{5}}
\begin{bmatrix}
-1 & 1 & 1 & 1 & 1\\
1 & -1 & 1 & 1 & 1\\
1 & 1 &- 1 & 1 & 1\\
1 & 1 & 1 & -1 & 1 \\
1 & 1 & 1 & 1 & -1\\
\end{bmatrix}\]

\subsection{$M=5$ vectors in $\R^5$ at angle $1/5$ 
and is circulant}

\[
 \frac{1}{\sqrt{5}}
\begin{bmatrix}
1 & 1 & 1 & -1 & 1\\
1 & 1 & 1 & 1 & -1\\
-1 & 1 & 1 & 1 & 1\\
1 & -1 & 1 & 1 & 1 \\
1 & 1 & -1 & 1 & 1
\end{bmatrix}\]

\subsection{$M=4$ vectors in $\R^7$ at angle $1/5$}
\[ \begin{bmatrix}
1 & 2 & 3 & 4 & 5 & 6 & 7 \\
\sqrt{\frac{4}{10}} & \sqrt{\frac{1}{10}} &
\sqrt{\frac{1}{10}}&\sqrt{\frac{1}{10}}&\sqrt{\frac{1}{10}}&
\sqrt{\frac{1}{10}}&\sqrt{\frac{1}{10}}\\
\sqrt{\frac{4}{10}} & -\sqrt{\frac{1}{10}} &
-\sqrt{\frac{1}{10}}&-\sqrt{\frac{1}{10}}&-\sqrt{\frac{1}{10}}&
\sqrt{\frac{1}{10}}&\sqrt{\frac{1}{10}}\\
\sqrt{\frac{4}{10}} & \sqrt{\frac{1}{10}} &
-\sqrt{\frac{1}{10}}&\sqrt{\frac{1}{10}}&-\sqrt{\frac{1}{10}}&
-\sqrt{\frac{1}{10}}&-\sqrt{\frac{1}{10}}\\
\sqrt{\frac{4}{10}} & -\sqrt{\frac{1}{10}} &
-\sqrt{\frac{1}{10}}&\sqrt{\frac{1}{10}}&\sqrt{\frac{1}{10}}&
-\sqrt{\frac{1}{10}}&-\sqrt{\frac{1}{10}} 
\end{bmatrix}\]

\subsection{$M=8$ vectors in $R^7$ at angle $1/5$}

\[\sqrt{\frac{1}{10}} \begin{bmatrix}
1&1&1&1& \sqrt{2}&\sqrt{2}&\sqrt{2}\\
1&1&-1&-1& \sqrt{2}&-\sqrt{2}&\sqrt{2}\\
1&-1&1&-1& \sqrt{2}&-\sqrt{2}&-\sqrt{2}\\
1&-1&-1&1& -\sqrt{2}&-\sqrt{2}&\sqrt{2}\\
1&1&1&1& -\sqrt{2}&-\sqrt{2}&-\sqrt{2}\\
1&1&-1&-1& -\sqrt{2}&\sqrt{2}&-\sqrt{2}\\
1&-1&1&-1& -\sqrt{2}&\sqrt{2}&\sqrt{2}\\
1&-1&-1&1& \sqrt{2}&\sqrt{2}&-\sqrt{2}\\
\end{bmatrix}\]

\subsection{$M=8$ vectors in $\R^8$ at angle $1/5$}

\[ 
\sqrt{\frac{1}{10}}\begin{bmatrix}
1&1&1&1&1&1&1&\sqrt{3}\\
1&1&-1&-1&-1&-1&1&\sqrt{3}\\ 
1&-1&1&-1&1&-1&-1&\sqrt{3}\\
1&-1&-1&1&-1&1&-1&\sqrt{3}\\
1&1&1&1&-1&-1&-1&-\sqrt{3}\\
1&1&-1&-1&1&1&-1&-\sqrt{3}\\
1&-1&1&-1&-1&1&1&-\sqrt{3}\\
1&-1&-1&1&1&-1&1&-\sqrt{3}\\
\end{bmatrix}\] 

\pagebreak

\noindent {\bf (B)  Specific Dimensions}
\vskip12pt

\subsection{$M=28$ vectors in $\R^{14}$ at angle $1/5$}

This example is built up from building blocks.
The first is:

\[ BB2= \begin{bmatrix}
\sqrt{\frac{2}{5}}& \sqrt{\frac{1}{5}}& \sqrt{\frac{1}{5}}
& \sqrt{\frac{1}{5}}\\
\sqrt{\frac{2}{5}}& -\sqrt{\frac{1}{5}}& -\sqrt{\frac{1}{5}}
& \sqrt{\frac{1}{5}}\\
\sqrt{\frac{2}{5}}& -\sqrt{\frac{1}{5}}& \sqrt{\frac{1}{5}}
& -\sqrt{\frac{1}{5}}\\
\sqrt{\frac{2}{5}}& \sqrt{\frac{1}{5}}& -\sqrt{\frac{1}{5}}
& -\sqrt{\frac{1}{5}}
\end{bmatrix}\]

This matrix will be spread out and the columns represented
by "bullets".

The second is:

\[ BB3=\begin{bmatrix}
1 & 2 & 3 & 4 & 5 & 6 & 7 \\
\sqrt{\frac{4}{10}} & \sqrt{\frac{1}{10}} &
\sqrt{\frac{1}{10}}&\sqrt{\frac{1}{10}}&\sqrt{\frac{1}{10}}&
\sqrt{\frac{1}{10}}&\sqrt{\frac{1}{10}}\\
\sqrt{\frac{4}{10}} & -\sqrt{\frac{1}{10}} &
-\sqrt{\frac{1}{10}}&-\sqrt{\frac{1}{10}}&-\sqrt{\frac{1}{10}}&
\sqrt{\frac{1}{10}}&\sqrt{\frac{1}{10}}\\
\sqrt{\frac{4}{10}} & \sqrt{\frac{1}{10}} &
\sqrt{\frac{1}{10}}&-\sqrt{\frac{1}{10}}&-\sqrt{\frac{1}{10}}&
-\sqrt{\frac{1}{10}}&-\sqrt{\frac{1}{10}}\\
\sqrt{\frac{4}{10}} & -\sqrt{\frac{1}{10}} &
-\sqrt{\frac{1}{10}}&\sqrt{\frac{1}{10}}&\sqrt{\frac{1}{10}}&
-\sqrt{\frac{1}{10}}&-\sqrt{\frac{1}{10}} 
\end{bmatrix}\]

For this matrix we will just put the first row in.

Now we use our dot trick to piece these together.
\begin{equation*}\setcounter{MaxMatrixCols}{14}
\begin{aligned}[b]\begin{array}{c}
BB3\backslash 4 \\[8pt]
BB2\backslash 4\\[8pt]
BB2 \backslash 4   \\[8pt]
BB2 \backslash 4  \\[8pt]
BB2 \backslash 4   \\[8pt]
BB2 \backslash 4   \\[8pt]
BB2 \backslash 4  \end{array}\end{aligned}
\begin{aligned}[b]&\quad\begin{array}{p{.65cm}p{.65cm}p{.65cm}p{.65cm}p{.65cm}p{.65cm}p{.65cm}p{.25cm}p{.25cm}p{.25cm}p{.25cm}p{.25cm}p{.25cm}p{.25cm}}
1 & 2 & 3 & 4 & 5 & 6 & 7 & 8 & 9 & 10 & 11 & 12 & 13 & 14\end{array}\\
&\left[\begin{array}{p{.65cm}p{.65cm}p{.65cm}p{.65cm}p{.65cm}p{.65cm}p{.65cm}p{.25cm}p{.25cm}p{.25cm}p{.25cm}p{.25cm}p{.25cm}p{.25cm}}
 $ \sqrt{\frac{4}{10}}$ & $\sqrt{\frac{1}{10}}$ & $\sqrt{\frac{1}{10}}
$&$\sqrt{\frac{1}{10}}$&$\sqrt{\frac{1}{10}}$& $\sqrt{\frac{1}{10}}
$&$\sqrt{\frac{1}{10}}$  & & & & & &  \\
  & $\sqrt{\frac{2}{5}}$ & & & & & & $\bullet$ & $\bullet $& $\bullet 
$& & & & \\
 & & $\sqrt{\frac{2}{5}}$ & & & & & & & $\bullet$ & $\bullet$ & $\bullet$ & & \\
 & & & $\sqrt{\frac{2}{5}}$& & & & $\bullet$ & & & & $\bullet$ & $\bullet$&\\
 & & & & $\sqrt{\frac{2}{5}}$ & & & &$\bullet$ & & $\bullet$& & $\bullet$& \\
 & & & & &$\sqrt{\frac{2}{5}}$ & & & & $\bullet$ & & & $\bullet$ & $\bullet $\\
  & & & & & &$\sqrt{\frac{2}{5}}$& $\bullet$& & & $\bullet$& & & $\bullet$
\end{array}\right]\end{aligned}\end{equation*}

\pagebreak

\subsection{$M=30$ in $\R^{15}$ at angle $1/5$}

We use a Building Block and put "bullets" each place it
occurs.

\[ BB4=
 \frac{1}{\sqrt{5}}
\begin{bmatrix}
-1 & 1 & 1 & 1 & 1\\
1 & -1 & 1 & 1 & 1\\
1 & 1 &- 1 & 1 & 1\\
1 & 1 & 1 & -1 & 1 \\
1 & 1 & 1 & 1 & -1\\
\end{bmatrix}\]

\begin{equation*}\begin{aligned}[b]&\quad
& \left[ \begin{array}{cccccccccccccccc}
&1&2&3&4&5&6&7&8&9&10&11&12&13&14&15\\
BB4\backslash 4&\bullet & \bullet & \bullet & \bullet & \bullet 
&&&&&&&&&&\\
BB4\backslash 4&&&&&\bullet &\bullet& \bullet & \bullet & \bullet 
&&&&&&\\
BB4\backslash 4&\bullet &&&&&&&&\bullet & \bullet & \bullet & \bullet 
&&&\\
BB4\backslash 4&& \bullet &&&&\bullet &&&& \bullet &&& 
\bullet &\bullet  \\
BB4\backslash 4&&& \bullet &&&&\bullet &&&&\bullet &&\bullet &&\bullet\\
BB4\backslash 4&&&& \bullet &&&&\bullet &&&&\bullet &&\bullet & \bullet \\
\end{array}\right] \end{aligned}
\end{equation*} 

\pagebreak

\subsection{$M=36$ vectors in $\R^{15}$ at angle $1/5$}

First we reverse one of our Building Blocks:

\[ BB5= \begin{bmatrix}
\sqrt{\frac{1}{5}}& \sqrt{\frac{1}{5}}& \sqrt{\frac{1}{5}}
& \sqrt{\frac{2}{5}}\\
\sqrt{\frac{1}{5}}& -\sqrt{\frac{1}{5}}& -\sqrt{\frac{1}{5}}
& \sqrt{\frac{2}{5}}\\
\sqrt{\frac{1}{5}}& -\sqrt{\frac{1}{5}}& \sqrt{\frac{1}{5}}
& -\sqrt{\frac{2}{5}}\\
\sqrt{\frac{1}{5}}& \sqrt{\frac{1}{5}}& -\sqrt{\frac{1}{5}}
& -\sqrt{\frac{2}{5}}
\end{bmatrix}\]

\[ BB6=
\sqrt{\frac{1}{10}}\begin{bmatrix}
1&1&1&1&1&1&1&\sqrt{3}\\
1&1&-1&-1&-1&-1&1&\sqrt{3}\\ 
1&-1&1&-1&1&-1&-1&\sqrt{3}\\
1&-1&-1&1&-1&1&-1&\sqrt{3}\\
1&1&1&1&-1&-1&-1&-\sqrt{3}\\
1&1&-1&-1&1&1&-1&-\sqrt{3}\\
1&-1&1&-1&-1&1&1&-\sqrt{3}\\
1&-1&-1&1&1&-1&1&-\sqrt{3}\\
\end{bmatrix}\]

\begin{equation*}
\begin{aligned}[b]&\quad
\left[ 
 \begin{array}{cccccccccccccccc}
&1&2&3&4&5&6&7&8&9&10&11&12&13&14&15\\
BB5 \ 4&\bullet&\bullet&\bullet & & & & &\sqrt{\frac{2}{5}}
&&&&&&&\\
BB5\ 4&&&\bullet &\bullet & \bullet& & & &\sqrt{\frac{2}{5}} &&&&&&\\
BB5\ 4&\bullet &&&&\bullet&\bullet &&&&\sqrt{\frac{2}{5}}&&&&&\\
BB5\ 4&& \bullet && \bullet && \bullet &&&&&\sqrt{\frac{2}{5}}&&&&\\
BB5\ 4&\bullet & & &\bullet && &\bullet &&&&&\sqrt{\frac{2}{3}}&&&\\
BB5\ 4&& \bullet &&& \bullet &&& \bullet &&&&&\sqrt{\frac{2}{5}} &&\\
BB5\ 4&& & \bullet & &&\bullet & \bullet &&&&&&&\sqrt{\frac{2}{5}} & \\
BB6\ 8&
&&&&&&&\sqrt{\frac{1}{10}}& \sqrt{\frac{1}{10}}&\sqrt{\frac{1}{10}}&
\sqrt{\frac{1}{10}}&\sqrt{\frac{1}{10}}&\sqrt{\frac{1}{10}}&
\sqrt{\frac{1}{10}}&\sqrt{\frac{3}{10}}\\
\end{array}\right] \end{aligned}
\end{equation*}

\pagebreak

\subsection{$M=40$ vectors in $\R^{16}$ at angle $1/5$}

For this we need a new building block.

\[ BB9=
\begin{bmatrix}
\sqrt{\frac{1}{5}}&
-\sqrt{\frac{1}{10}} & -\sqrt{\frac{1}{10}}&-\sqrt{\frac{1}{10}}&
\sqrt{\frac{1}{10}}&\sqrt{\frac{2}{10}}\\
\sqrt{\frac{1}{5}}&
\sqrt{\frac{1}{10}} & \sqrt{\frac{1}{10}}&\sqrt{\frac{1}{10}}&
-\sqrt{\frac{1}{10}}&\sqrt{\frac{2}{10}}\\
-\sqrt{\frac{1}{5}}&
-\sqrt{\frac{1}{10}} & \sqrt{\frac{1}{10}}&-\sqrt{\frac{1}{10}}&
-\sqrt{\frac{1}{10}}&\sqrt{\frac{2}{10}}\\
-\sqrt{\frac{1}{5}}&
\sqrt{\frac{1}{10}} & -\sqrt{\frac{1}{10}}&\sqrt{\frac{1}{10}}&
\sqrt{\frac{1}{10}}&\sqrt{\frac{2}{10}}\\
\end{bmatrix}\]

In the next chart, the first seven rows are building block 5
and contain four vectors each.  Row 8 is Building block 6 and
has eight vectors.  And Row 9 is Building block 9 and has
four vectors. 
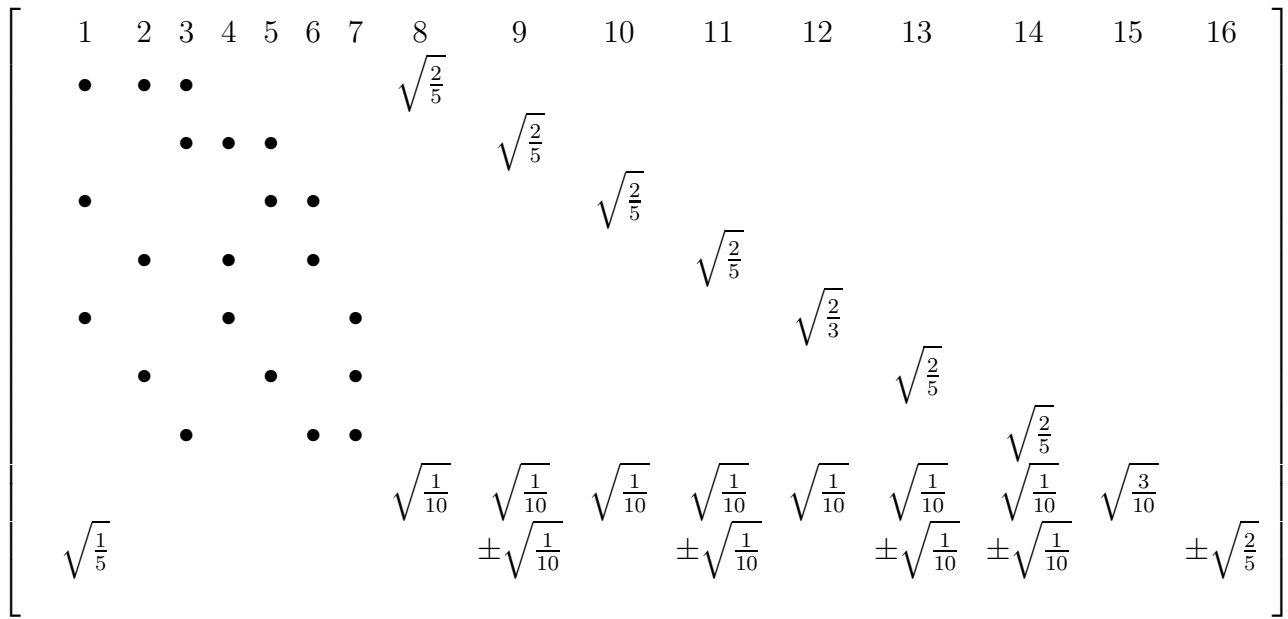
\begin{equation*}
\begin{aligned}[b]&\quad
\left[ 
 \scriptscriptstyle\begin{array}{ccccccccccccccccc}
&1&2&3&4&5&6&7&8&9&10&11&12&13&14&15&16\\
&\bullet&\bullet&\bullet & & & & &\sqrt{\frac{2}{5}}
&&&&&&&&\\
&&&\bullet &\bullet & \bullet& & & &\sqrt{\frac{2}{5}} &&&&&&&\\
&\bullet &&&&\bullet&\bullet &&&&\sqrt{\frac{2}{5}}&&&&&&\\
&& \bullet && \bullet && \bullet &&&&&\sqrt{\frac{2}{5}}&&&&&\\
&\bullet & & &\bullet && &\bullet &&&&&\sqrt{\frac{2}{3}}&&&&\\
&& \bullet &&& \bullet &&\bullet &&&&&&\sqrt{\frac{2}{5}} &&&\\
&& & \bullet & &&\bullet & \bullet &&&&&&&\sqrt{\frac{2}{5}} && \\
&
&&&&&&&\sqrt{\frac{1}{10}}& \sqrt{\frac{1}{10}}&\sqrt{\frac{1}{10}}&
\sqrt{\frac{1}{10}}&\sqrt{\frac{1}{10}}&\sqrt{\frac{1}{10}}&
\sqrt{\frac{1}{10}}&\sqrt{\frac{3}{10}}&\\
&\sqrt{\frac{1}{5}}&&&&&&&&\pm \sqrt{\frac{1}{10}}
&&\pm \sqrt{\frac{1}{10}}&
&\pm \sqrt{\frac{1}{10}}&\pm \sqrt{\frac{1}{10}}&&\pm \sqrt{\frac{2}{5}}\\
\scriptscriptstyle\end{array}\right] \end{aligned}
\end{equation*}

On the next page we will give a hint about how to prove this
is an equiangular line set.

\pagebreak

The first 36 vectors represent 36 equiangular lines in $\R^{15}$ and
this is easily checked by sight.  Also, the last four vectors have
$\sqrt{\frac{1}{5}}$ in the first position and this hits exactly
the groups in rows 1,3 and 5 in eactly one position - the first position
and this yields $1/5$.  The last group hits the elements in rows 2,4,6, and
7 in exactly one place and produces
\[ \sqrt{\frac{2}{5}} \cdot \sqrt{\frac{1}{10}} = \frac{1}{5}.\]

So we only need to check how the last four vectors interact with the
eight just before them.  But these vectors all have exactly four
coordinates in common (9,11,13, and 14) and they are respectively:

\[ \begin{bmatrix}
+&+&+&+\\
+&-&-&+\\
-&-&-&-\\
-&+&+&-\\
+&+&-&-\\
+&-&+&-\\
-&-&+&+\\
-&+&-&+\\
\end{bmatrix} \]

and

\[ \begin{bmatrix}
-&-&-&+\\
+&+&+&-\\
-&+&-&-\\
+&-&+&+\\
\end{bmatrix} \]

Now, a visual check shows that the inner products of the rows in
the second group with the rows in the first group will yield
\[
\pm \left [ \sqrt{\frac{1}{10}}\cdot \sqrt{\frac{1}{10}}
+ \sqrt{\frac{1}{10}}\cdot \sqrt{\frac{1}{10}}\right ] =
\pm \frac{1}{5}.\]

\pagebreak

\subsection{$M=60$ vectors in $\R^{19}$ at angle $1/5$}

We use a Building Block and put "bullets" each place it
occurs.

\[ BB4=
 \frac{1}{\sqrt{5}}
\begin{bmatrix}
-1 & 1 & 1 & 1 & 1\\
1 & -1 & 1 & 1 & 1\\
1 & 1 &- 1 & 1 & 1\\
1 & 1 & 1 & -1 & 1 \\
1 & 1 & 1 & 1 & -1\\
\end{bmatrix}\]

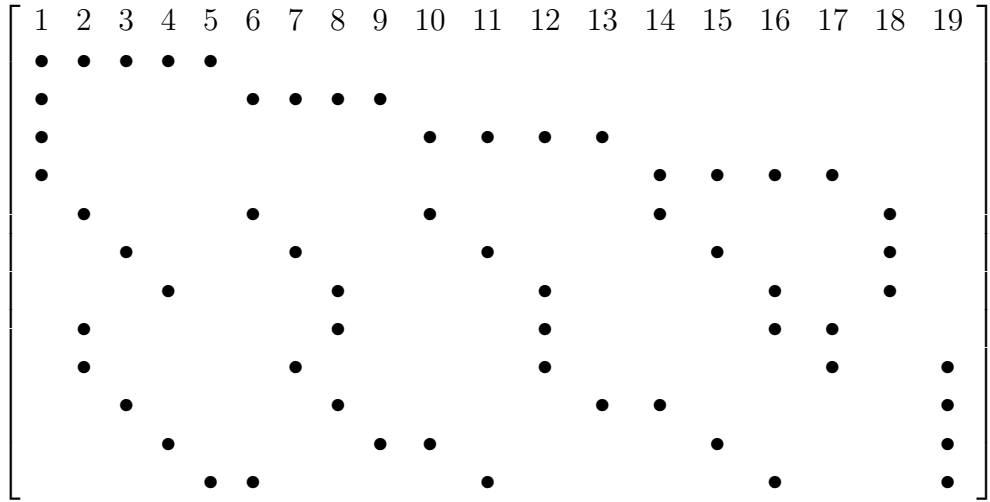
\begin{equation*}\begin{aligned}[b]&\quad
& \left[ \begin{array}{ccccccccccccccccccc}
1&2&3&4&5&6&7&8&9&10&11&12&13&14&15&16&17&18&19\\
\bullet & \bullet & \bullet & \bullet & \bullet 
&&&&&&&&&&&&&&\\
\bullet &&&&&\bullet& \bullet & \bullet & \bullet 
&&&&&&&&&&\\
\bullet &&&&&&&&&\bullet & \bullet & \bullet & \bullet 
&&&&&&\\
\bullet &&&&&&&&&&&&&\bullet & \bullet & \bullet & \bullet 
&&\\
& \bullet &&&&\bullet &&&& \bullet &&&& 
\bullet &&&&\bullet & \\
&& \bullet &&&&\bullet &&&&\bullet &&&&\bullet &&&\bullet&\\
&&& \bullet &&&&\bullet &&&&\bullet &&&&\bullet && \bullet &\\
& \bullet &&&&&&\bullet &&&&\bullet &&&& \bullet & \bullet &\\
& \bullet &&&&&\bullet &&&&& \bullet &&&&& 
\bullet &&\bullet \\
&& \bullet &&&&&\bullet &&&&&\bullet &\bullet &&&&&\bullet\\
&&& \bullet &&&&&\bullet &\bullet &&&&&\bullet &&&& \bullet \\
&&&& \bullet &\bullet &&&&&\bullet &&&&& \bullet &&& \bullet \\
\end{array}\right] \end{aligned}
\end{equation*}

\pagebreak

\subsection{$M=45$ vectors in $\R^{21}$ at angle $1/5$}

We use two building blocks.  The first is a $4\times 4$
orthogonal matrix with an extra column - and it is
represented in our matrix by "bullets".

\[ BB7=\sqrt{\frac{1}{5}}\begin{bmatrix}
1&1&1&1&1\\
1&1&1&-1&-1\\
1&1&-1&1&-1\\
1&1&-1&-1&1
\end{bmatrix}\]

The second building block is a $5\times 5$ matrix - and
it is represented by "*".

\[ BB8=
 \frac{1}{\sqrt{5}}
\begin{bmatrix}
1 & 1 & 1 & 1 & 1\\
-1 & -1 & -1 & 1 & 1\\
-1 & 1 & 1 & -1 & 1\\
-1 & -1 & 1 & 1 & -1 \\
1 & -1 & 1 & -1 & 1\\
\end{bmatrix}\]

\begin{equation*}\begin{aligned}[b]&\quad
& \left[ \begin{array}{ccccccccccccccccccccc}
1&2&3&4&5&6&7&8&9&10&11&12&13&14&15&16&17&18&19&20&21\\
\bullet & \bullet & \bullet & \bullet & \bullet 
&&&&&&&&&&&&&&&&\\
\bullet &&&&&\bullet& \bullet & \bullet & \bullet 
&&&&&&&&&&&&\\
\bullet &&&&&&&&&\bullet & \bullet & \bullet & \bullet 
&&&&&&&&\\
\bullet &&&&&&&&&&&&&\bullet & \bullet & \bullet & \bullet 
&&&&\\
\bullet &&&&&&&&&&&&&&&&&\bullet & \bullet & 
\bullet & \bullet\\
& * &&&& * &&&& * &&&& * &&&& * &&&\\
& * &&&&& * &&&& * &&&& * &&&& * &&\\
& * &&&&&& * &&&& * &&&& * &&&& * &\\
& * &&&&&&& * &&&& * &&&& * &&&& *\\
\end{array}\right] \end{aligned}
\end{equation*} 

\pagebreak

\subsection{$M=45$ vectors in $\R^{21}$ at angle $1/5$}

We use the a Building Block and put "bullets" each
place it appears.

\[ BB4=
 \frac{1}{\sqrt{5}}
\begin{bmatrix}
-1 & 1 & 1 & 1 & 1\\
1 & -1 & 1 & 1 & 1\\
1 & 1 &- 1 & 1 & 1\\
1 & 1 & 1 & -1 & 1 \\
1 & 1 & 1 & 1 & -1\\
\end{bmatrix}\]

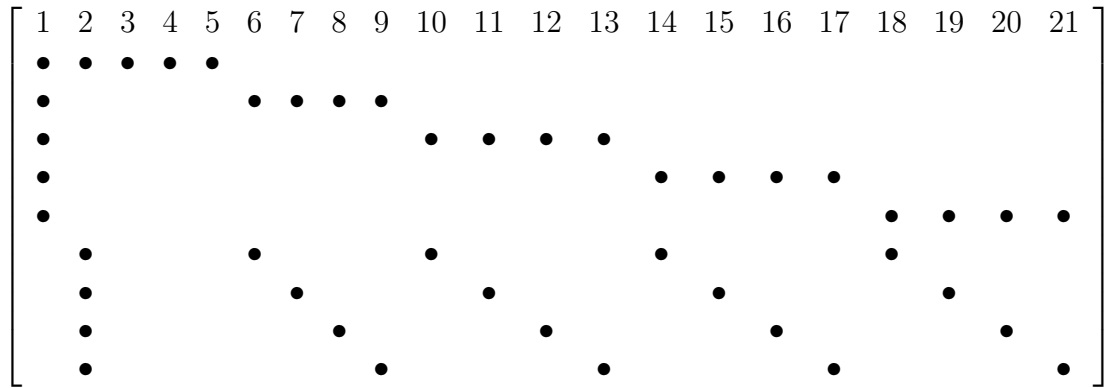
\begin{equation*}\begin{aligned}[b]&\quad
&\left[ \begin{array}{ccccccccccccccccccccc}
1&2&3&4&5&6&7&8&9&10&11&12&13&14&15&16&17&18&19&20&21\\
\bullet & \bullet & \bullet & \bullet & \bullet 
&&&&&&&&&&&&&&&&\\
\bullet &&&&&\bullet& \bullet & \bullet & \bullet 
&&&&&&&&&&&&\\
\bullet &&&&&&&&&\bullet & \bullet & \bullet & \bullet 
&&&&&&&&\\
\bullet &&&&&&&&&&&&&\bullet & \bullet & \bullet & \bullet 
&&&&\\
\bullet &&&&&&&&&&&&&&&&&\bullet & \bullet & 
\bullet & \bullet\\
& \bullet &&&&\bullet &&&&\bullet &&&&\bullet &&&&\bullet&&&\\
& \bullet &&&&&\bullet &&&&\bullet &&&&\bullet &&&& \bullet &&\\
& \bullet &&&&&&\bullet &&&&\bullet &&&& \bullet &&&& \bullet &\\
& \bullet &&&&&&&\bullet &&&&\bullet &&&& \bullet &&&& \bullet\\
\end{array}\right] \end{aligned}
\end{equation*} 

\pagebreak

\noindent {\bf (C)  General Sets}
\vskip12pt

\subsection{$N-1$ vectors in $\R^N$ at angle $1/5$}

\[ \begin{bmatrix}
1&\sqrt{\frac{1}{5}}& \sqrt{\frac{2}{5}}& \sqrt{\frac{2}{5}}&&&&\cdots\\
2&\sqrt{\frac{1}{5}}& \sqrt{\frac{2}{5}}& -\sqrt{\frac{2}{5}}&&&&\cdots\\ 
3&\sqrt{\frac{1}{5}}&& \sqrt{\frac{2}{5}}& \sqrt{\frac{2}{5}}&&&\cdots\\
4&\sqrt{\frac{1}{5}}&& \sqrt{\frac{2}{5}}& -\sqrt{\frac{2}{5}}&&&\cdots\\
5&\sqrt{\frac{1}{5}}&&& \sqrt{\frac{2}{5}}& \sqrt{\frac{2}{5}}&&\cdots\\
6&\sqrt{\frac{1}{5}}&&& \sqrt{\frac{2}{5}}& -\sqrt{\frac{2}{5}}&&\cdots\\
\vdots& \vdots & \vdots & \vdots & \vdots & \vdots & \vdots & \vdots 
\end{bmatrix} \]

\subsection{$N-1$ vectors in $\R^N$ at angle $1/5$}

\[ \begin{bmatrix}
1&\sqrt{\frac{1}{5}}& \sqrt{\frac{1}{5}}& \sqrt{\frac{3}{5}}&&&&\cdots\\
2&\sqrt{\frac{1}{5}}& \sqrt{\frac{1}{5}}& -\sqrt{\frac{3}{5}}&&&&\cdots\\ 
3&\sqrt{\frac{1}{5}}&& \sqrt{\frac{1}{5}}& \sqrt{\frac{3}{5}}&&&\cdots\\
4&\sqrt{\frac{1}{5}}&& \sqrt{\frac{1}{5}}& -\sqrt{\frac{3}{5}}&&&\cdots\\
5&\sqrt{\frac{1}{5}}&&& \sqrt{\frac{1}{5}}& \sqrt{\frac{3}{5}}&&\cdots\\
6&\sqrt{\frac{1}{5}}&&& \sqrt{\frac{1}{5}}& -\sqrt{\frac{3}{5}}&&\cdots\\
\vdots& \vdots & \vdots & \vdots & \vdots & \vdots & \vdots & \vdots 
\end{bmatrix} \]

\pagebreak

\subsection{$M=4N$ vectors in $\R^{3N+1}$ at angle $1/5$}

We use our building block 

\[ BB5= \begin{bmatrix}
\sqrt{\frac{1}{5}}& \sqrt{\frac{1}{5}}& \sqrt{\frac{1}{5}}
& \sqrt{\frac{2}{5}}\\
\sqrt{\frac{1}{5}}& -\sqrt{\frac{1}{5}}& -\sqrt{\frac{1}{5}}
& \sqrt{\frac{2}{5}}\\
\sqrt{\frac{1}{5}}& -\sqrt{\frac{1}{5}}& \sqrt{\frac{1}{5}}
& -\sqrt{\frac{2}{5}}\\
\sqrt{\frac{1}{5}}& \sqrt{\frac{1}{5}}& -\sqrt{\frac{1}{5}}
& -\sqrt{\frac{2}{5}}
\end{bmatrix}\]

We use "bullets" to indicate where these columns go.

\[ \begin{bmatrix}
\bullet & \bullet & \bullet & \bullet & &&&&&&\cdots \\
\bullet &&&&\bullet & \bullet & \bullet &&&&\cdots \\
\bullet &&&&&&& \bullet & \bullet & \bullet & \cdots\\
\vdots & \vdots & \vdots & \vdots & \vdots & 
\vdots & \vdots & \vdots & \vdots & \vdots &\cdots
\end{bmatrix} \]

\pagebreak
\section{Angle $1/7$}

\noindent {\bf (A)  Building Blocks}
\vskip12pt

\subsection{$M=7$ vectors in $\R^7$ at angle $1/7$}

\begin{equation*}\frac{1}{\sqrt{7}}
\begin{aligned}
&\left[ \begin{array}{ccccccc}
1&1&1&1&1&1&1\\
-1&-1&-1&1&1&1&1\\
1&1&-1&-1&-1&1&1\\
1&-1&1&1&-1&-1&1\\
1&1&-1&1&1&-1&-1\\
-1&1&1&1&1&-1&1\\
-1&1&-1&1&-1&-1&1\\
\end{array}\right] \end{aligned}\]

\subsection{$M=8$ vectors in $\R^7$ at angle $1/7$}
\[\sqrt{\frac{1}{10}} \begin{bmatrix}
1&1&1&1& \sqrt{2}&\sqrt{2}&\sqrt{2}\\
1&1&-1&-1& \sqrt{2}&-\sqrt{2}&\sqrt{2}\\
1&-1&1&-1& \sqrt{2}&-\sqrt{2}&-\sqrt{2}\\
1&-1&-1&1& -\sqrt{2}&-\sqrt{2}&\sqrt{2}\\
1&1&1&1& -\sqrt{2}&-\sqrt{2}&-\sqrt{2}\\
1&1&-1&-1& -\sqrt{2}&\sqrt{2}&-\sqrt{2}\\
1&-1&1&-1& -\sqrt{2}&\sqrt{2}&\sqrt{2}\\
1&-1&-1&1& \sqrt{2}&\sqrt{2}&-\sqrt{2}\\
\end{bmatrix}\]
\pagebreak

\section{$M=2N$ vectors in $\R^N$}

\subsection{$M=6$ vectors in $\R^3$ at angle $1/\sqrt{5}$}

\[
\begin{bmatrix}
0 & \sqrt{\frac{5-\sqrt{5}}{10}} & \sqrt{\frac{5+\sqrt{5}}{10}}\\
0 & -\sqrt{\frac{5-\sqrt{5}}{10}} & \sqrt{\frac{5+\sqrt{5}}{10}}\\
\sqrt{\frac{5-\sqrt{5}}{10}} & \sqrt{\frac{5+\sqrt{5}}{10}} & 0\\
-\sqrt{\frac{5-\sqrt{5}}{10}} & \sqrt{\frac{5+\sqrt{5}}{10}} & 0\\
\sqrt{\frac{5+\sqrt{5}}{10}} & 0 & \sqrt{\frac{5-\sqrt{5}}{10}}\\
\sqrt{\frac{5+\sqrt{5}}{10}} & 0 & -\sqrt{\frac{5-\sqrt{5}}{10}}
\end{bmatrix}\]

\subsection{$M=10$ vectors in $\R^5$ at angle $1/3$}

\[\begin{bmatrix}
\sqrt{\frac{1}{3}} & \sqrt{\frac{2}{3}} & 0 & 0 & 0\\
-\sqrt{\frac{1}{3}} & \sqrt{\frac{2}{3}} & 0 & 0 & 0\\
\sqrt{\frac{1}{3}} & 0 & \sqrt{\frac{2}{3}} & 0 & 0\\
- \sqrt{\frac{1}{3}} & 0 &\sqrt{\frac{2}{3}} & 0 & 0\\
\sqrt{\frac{1}{3}} & 0 & 0 & \sqrt{\frac{2}{3}} & 0\\
-\sqrt{\frac{1}{3}} & 0 & 0 &\sqrt{\frac{2}{3}} & 0\\
0 & \frac{1}{3}\sqrt{\frac{3}{2}} & \frac{1}{3} \sqrt{\frac{3}{2}} &
\frac{1}{3}\sqrt{\frac{3}{2}} & \sqrt{\frac{1}{2}}\\
0 & \frac{1}{3}\sqrt{\frac{3}{2}} & -\frac{1}{3} \sqrt{\frac{3}{2}} &
\frac{1}{3}\sqrt{\frac{3}{2}} & -\sqrt{\frac{1}{2}}\\
0 & -\frac{1}{3}\sqrt{\frac{3}{2}} & -\frac{1}{3} \sqrt{\frac{3}{2}} &
\frac{1}{3}\sqrt{\frac{3}{2}} & \sqrt{\frac{1}{2}}\\
0 & \frac{1}{3}\sqrt{\frac{3}{2}} & -\frac{1}{3} \sqrt{\frac{3}{2}} &
-\frac{1}{3}\sqrt{\frac{3}{2}} & \sqrt{\frac{1}{2}}\\
\end{bmatrix}\]

\subsection{$M=12$ vectors in $\R^6$ at angle $1/3$}

\begin{equation*}
\begin{aligned}[b] &\quad 
\left[ \begin{array}{cccccc}
1&2&3&4&5&6\\
\sqrt{\frac{1}{3}}& \sqrt{\frac{2}{3}}&&&&\\
\sqrt{\frac{1}{3}}&-\sqrt{\frac{2}{3}}&&&&\\
\sqrt{\frac{1}{3}} & & \sqrt{\frac{2}{3}}&&&\\
\sqrt{\frac{1}{3}}& & -\sqrt{\frac{2}{3}}&&&\\
\sqrt{\frac{1}{3}}& & & \sqrt{\frac{2}{3}}&&\\
\sqrt{\frac{1}{3}}& & & -\sqrt{\frac{2}{3}}&&\\
\sqrt{\frac{1}{3}}& & & & \sqrt{\frac{2}{3}}&\\
\sqrt{\frac{1}{3}}& & & &-\sqrt{\frac{2}{3}}&\\
& \frac{1}{\sqrt{6}}& \frac{1}{\sqrt{6}}&\frac{1}{\sqrt{6}}
&\frac{1}{\sqrt{6}}&\frac{1}{\sqrt{3}}\\
&\frac{1}{\sqrt{6}}& \frac{1}{\sqrt{6}}&\frac{-1}{\sqrt{6}}
&\frac{-1}{\sqrt{6}}&\frac{1}{\sqrt{3}}\\
&\frac{1}{\sqrt{6}}& \frac{-1}{\sqrt{6}}&\frac{1}{\sqrt{6}}
&\frac{-1}{\sqrt{6}}&\frac{1}{\sqrt{3}}\\
&\frac{1}{\sqrt{6}}& \frac{-1}{\sqrt{6}}&\frac{-1}{\sqrt{6}}
&\frac{1}{\sqrt{6}}&\frac{1}{\sqrt{3}}\\
\end{array}\right] \end{aligned}\]

\subsection{$M=28$ vectors in $\R^{14}$ at angle $1/3$}

First, we use our {\em building block} for $4$ vectors
in $\R^3$ at angle $1/3$ where $+$ means $1$ and $-$ means
$-1$.

\[ BB1 =\sqrt{\frac{1}{3}}\begin{bmatrix}
+ & + & +\\
- & + & + \\
+ & - & + \\
+ & + & - \\
\end{bmatrix}\]
Now, we put these groups of four vectors into
our chart by just putting a {\em dot} where the
columns go.

\begin{equation*}\begin{aligned}[b]\begin{array}{c}
BB1\  4 \\
BB1\  4\\
BB1\  4 \\
BB1\ 4\\
BB1\ 4 \\
BB1\ 4\\
BB1\ 4 \\
\end{array}\end{aligned}
\begin{aligned}[b] &\quad \begin{array}{ccccccc}
1 & 2 & 3 & 4 & 5 & 6 & 7\end{array}\\
& \left[ \begin{array}{ccccccc}
\bullet & \bullet & \bullet & &&&\\
& & \bullet & \bullet & \bullet & &\\
\bullet & & & &\bullet & \bullet &\\
& \bullet && \bullet && \bullet &\\
& & \bullet & & & \bullet & \bullet\\
\bullet & & & \bullet & & &\bullet\\
& \bullet & & & \bullet & & \bullet\\
\end{array}\right]\end{aligned}
\end{equation*}

\pagebreak

\subsection{$M=28$ vectors in $\R^{14}$ at angle $1/5$}

This example is built up from building blocks.
The first is:

\[ BB2= \begin{bmatrix}
\sqrt{\frac{2}{5}}& \sqrt{\frac{1}{5}}& \sqrt{\frac{1}{5}}
& \sqrt{\frac{1}{5}}\\
\sqrt{\frac{2}{5}}& -\sqrt{\frac{1}{5}}& -\sqrt{\frac{1}{5}}
& \sqrt{\frac{1}{5}}\\
\sqrt{\frac{2}{5}}& -\sqrt{\frac{1}{5}}& \sqrt{\frac{1}{5}}
& -\sqrt{\frac{1}{5}}\\
\sqrt{\frac{2}{5}}& \sqrt{\frac{1}{5}}& -\sqrt{\frac{1}{5}}
& -\sqrt{\frac{1}{5}}
\end{bmatrix}\]

This matrix will be spread out and the columns represented
by "bullets".

The second is:

\[ BB3=\begin{bmatrix}
1 & 2 & 3 & 4 & 5 & 6 & 7 \\
\sqrt{\frac{4}{10}} & \sqrt{\frac{1}{10}} &
\sqrt{\frac{1}{10}}&\sqrt{\frac{1}{10}}&\sqrt{\frac{1}{10}}&
\sqrt{\frac{1}{10}}&\sqrt{\frac{1}{10}}\\
\sqrt{\frac{4}{10}} & -\sqrt{\frac{1}{10}} &
-\sqrt{\frac{1}{10}}&-\sqrt{\frac{1}{10}}&-\sqrt{\frac{1}{10}}&
\sqrt{\frac{1}{10}}&\sqrt{\frac{1}{10}}\\
\sqrt{\frac{4}{10}} & \sqrt{\frac{1}{10}} &
-\sqrt{\frac{1}{10}}&\sqrt{\frac{1}{10}}&-\sqrt{\frac{1}{10}}&
-\sqrt{\frac{1}{10}}&-\sqrt{\frac{1}{10}}\\
\sqrt{\frac{4}{10}} & -\sqrt{\frac{1}{10}} &
-\sqrt{\frac{1}{10}}&\sqrt{\frac{1}{10}}&\sqrt{\frac{1}{10}}&
-\sqrt{\frac{1}{10}}&-\sqrt{\frac{1}{10}} 
\end{bmatrix}\]

For this matrix we will just put the first row in.

Now we use our dot trick to piece these together.
\begin{equation*}\setcounter{MaxMatrixCols}{14}
\begin{aligned}[b]\begin{array}{c}
BB3\backslash 4 \\[8pt]
BB2\backslash 4\\[8pt]
BB2 \backslash 4   \\[8pt]
BB2 \backslash 4  \\[8pt]
BB2 \backslash 4   \\[8pt]
BB2 \backslash 4   \\[8pt]
BB2 \backslash 4  \end{array}\end{aligned}
\begin{aligned}[b]&\quad\begin{array}{p{.65cm}p{.65cm}p{.65cm}p{.65cm}p{.65cm}p{.65cm}p{.65cm}p{.25cm}p{.25cm}p{.25cm}p{.25cm}p{.25cm}p{.25cm}p{.25cm}}
1 & 2 & 3 & 4 & 5 & 6 & 7 & 8 & 9 & 10 & 11 & 12 & 13 & 14\end{array}\\
&\left[\begin{array}{p{.65cm}p{.65cm}p{.65cm}p{.65cm}p{.65cm}p{.65cm}p{.65cm}p{.25cm}p{.25cm}p{.25cm}p{.25cm}p{.25cm}p{.25cm}p{.25cm}}
 $ \sqrt{\frac{4}{10}}$ & $\sqrt{\frac{1}{10}}$ & $\sqrt{\frac{1}{10}}$&$\sqrt{\frac{1}{10}}$&$\sqrt{\frac{1}{10}}$& $\sqrt{\frac{1}{10}}$&$\sqrt{\frac{1}{10}}$  & & & & & &  \\
  & $\sqrt{\frac{2}{5}}$ & & & & & & $\bullet$ & $\bullet $& $\bullet $& & & & \\
 & & $\sqrt{\frac{2}{5}}$ & & & & & & & $\bullet$ & $\bullet$ & $\bullet$ & & \\
 & & & $\sqrt{\frac{2}{5}}$& & & & $\bullet$ & & & & $\bullet$ & &\\
 & & & & $\sqrt{\frac{2}{5}}$ & & & &$\bullet$ & & $\bullet$& & $\bullet$& \\
 & & & & &$\sqrt{\frac{2}{5}}$ & & & & $\bullet$ & & & $\bullet$ & $\bullet $\\
  & & & & & &$\sqrt{\frac{2}{5}}$& $\bullet$& & & $\bullet$& & & $\bullet$\end{array}\right]\end{aligned}\end{equation*}

\pagebreak
\section{$M=N+1$ vectors in $\R^N$}

\begin{proposition}
For every natural number $N$, the $(N+1)\times (N+1)$ matrix below
represents $N+1$ unit norm vectors $\{f_m\}_{m=1}^{N+1}$ in
$\R^{N+1}$ giving an equiangular tight frame for an $N$-dimensional
space and satisfying:

(1)  $\|f_m\|=1$.

(2)  For any $n\not= m$, $\langle f_m,f_n\rangle = \frac{-1}{N}$.
\end{proposition}

\begin{proof}
For any $m$ we have
\[ \|f_m\|^2 = \frac{1}{N(N+1}[N^2+N] =1.\]
Also,
\[ \langle f_m,f_n\rangle = \frac{-2N+N-1}{N(N+1)}
= \frac{-(N+1)}{N(N+1)} = \frac{-1}{N}.\]
\end{proof}

\pagebreak

\section{Unitary Matrices}

\subsection{Circulant self-adjoint Matrix}

\[ \begin{bmatrix}
-b & a & a & \cdots & a\\
a & -b & a & \cdots & a\\
\vdots & \vdots & \vdots & \cdots & \vdots\\
a & a & a & \cdots & a\\
a & a & a & \cdots & -b
\end{bmatrix}
\]

\subsection{Circulant self-adjoint Unitary Matrix}
If we rotate the rows of the matrix we can get matrices 
for every $2N$-dimensional Hilbert
space yielding a circulant, self-adjoint unitary matrix.

\begin{equation*}\left[\begin{array}{cccc|cccc}
a & a & \dots & a & -b & a & \dots & a\\
a & a & \dots & a & a & -b & \dots &a\\
\vdots & \vdots & \vdots & \vdots & \vdots& \vdots & \vdots & \vdots\\
a & a & \dots & a & a & a & \dots & -b\\\hline
-b & a & \dots & a & a & a & \dots & a\\
a & -b & \dots & a & a & a &\dots & a\\
\vdots & \vdots & \vdots & \vdots & \vdots & \vdots & \vdots & \vdots\\
a & a & \dots & -b & a & a &\dots &a\end{array}\right]\end{equation*}

\begin{proof}
We have for any two rows $n\not= m$ of this matrix:
\[ \langle f_n,f_m\rangle = (2N-2)a2 - 2ab =
a[(2(N-1)a-2b].\]
So this is zero if $b=(N-1)a$.
\end{proof}

\subsection{Another example of circulant matrices}

\[ \begin{bmatrix}
a & -b & a & \cdots & a\\
a & a & -b & \cdots & a\\
\vdots & \vdots & \vdots & \cdots & \vdots\\
a & a & a & \cdots & -b\\
-b & a & a & \cdots & a
\end{bmatrix}
\]

\pagebreak

\section{Multiple angles}

\subsection{$M=2N$ vectors at two angles: $1/\sqrt{5},0$}
\begin{proposition}
For every $N$, $\R^N$ has $2N$ lines at three angles:
$\pm 1/\sqrt{5},0$.  These lines span
$\R^N$.
\end{proposition}

\begin{proof}
We let
\[ x = \sqrt{\frac{5-\sqrt{5}}{10}}\ \, y =
\sqrt{\frac{5+\sqrt{5}}{10}}.\]
We will just write down the vectors:
\[ \begin{bmatrix}
f_1 \\ f_2 \\ f_3 \\ f_4\\ f_5 \\ f_6 \\ \vdots \\ f_{2N-3}\\
f_{2N-2}\\ f_{2N-1}\\ f_{2N}.
\end{bmatrix}
=
\begin{bmatrix}
x & y & 0 & 0 &  \cdots & 0 & 0 & 0 \\
-x & y & 0 & 0 & \cdots & 0 & 0 & 0\\
0 & x & y & 0 & \cdots & 0 & 0 & 0 \\
0 & -x & y & 0 & \cdots & 0 & 0 & 0\\
0 & 0 & x & y & \cdots & 0 & 0 & 0\\
0 & 0 & -x & y & \cdots & 0 & 0 & 0\\
\vdots & \vdots & \vdots & \vdots & \vdots &
\vdots & \vdots & \vdots & \\
0 & 0 & 0 & 0 & \cdots & 0 & x & y \\
0 & 0 & 0 & 0 & \cdots & 0 & -x & y\\
y & 0 & 0 & 0 & \cdots & 0 & 0 & x\\
y & 0 & 0 & 0 & \cdots & 0 & 0 & -x
\end{bmatrix}
\]
It is obvious that these lines span $\R^N$.
\end{proof}

\begin{corollary}
For every $N$, the $2N$ unit vectors at angles
$1/\sqrt{5},0$ can be divided into two sets of
circulant vectors.
\end{corollary}

\begin{proof}
The two sets are:
\[ \begin{bmatrix}
f_1 \\ f_2 \\ f_3 \\ \vdots \\ f_N
\end{bmatrix}
=
\begin{bmatrix}
x & y & 0 & 0 & \cdots & 0 & 0\\
0 & x & y & 0 & \cdots & 0 & 0 \\
0 & 0 & x & y & \cdots & 0 & 0\\
\vdots & \vdots & \vdots & \vdots &
\cdots & 0 & 0 \\
0 & 0 & 0 & 0 & \cdots & x & y\\
y & 0 & 0 & 0 & \cdots & 0 & x
\end{bmatrix}\]

\[ \begin{bmatrix}
f_1 \\ f_2 \\ f_3 \\ \vdots \\ f_N
\end{bmatrix}
=
\begin{bmatrix}
-x & y & 0 & 0 & \cdots & 0 & 0\\
0 & -x & y & 0 & \cdots & 0 & 0 \\
0 & 0 & -x & y & \cdots & 0 & 0\\
\vdots & \vdots & \vdots & \vdots &
\cdots & 0 & 0 \\
0 & 0 & 0 & 0 & \cdots & -x & y\\
y & 0 & 0 & 0 & \cdots & 0 & -x
\end{bmatrix}\]
\end{proof}


\begin{thebibliography}{WW}
\bibitem{BK}
J. J. Benedetto and J. Kolesar,  
\emph{Geometric properties of Grassmannian frames for $R^2$ and $R^3$}, 
EURASIP J. Applied Signal Processing, (2006), 17 pages.

\bibitem{CRT}  P.G. Casazza, D. Redmond and J.C. Tremain,
\emph{Real equiangular tight frames}, in preparation.

\bibitem{CRT1}  P.G. Casazza, D. Redmond and J.C. Tremain,
\emph{Real equiangular frames}, Proceedings of CISS - Princeton
University (2008).

\bibitem{H}  J.~Haantjes, \emph{Equilateral point-sets in elliptic two-
and three-dimensional spaces}, Nieuw Arch. Wisk. {\bf 22} (1948)
355-362.


\bibitem{HP}  R. B.~Holmes and V. I.~Paulsen, \emph{Optimal frames for
erasures}, Linear Alg. and Applications {\bf 377} (204) 31-51.


\bibitem{LS}  P.~Lemmens and J.~Seidel, \emph{Equiangular lines},
Journal of Algebra {\bf 24} (1973) 494-512.

\bibitem{LiS}  J. H.~van Lint and J. J.~Seidel, \emph{Equiangular point
sets in elliptic geometry}, Proc. Nederl. Akad. Wetensch. Series A {\bf 69}
(1966) 335-348.



\bibitem{SH} T.~Strohmer and R. W.~Heath, \emph{Grassmannian frames
with applications to coding and communication}, Appl. Comp. Harmonic
Anal. {\bf 14} No. 3 (2003) 257-275.



\bibitem{W}  L. R.~Welch, \emph{Lower bounds on the maximum cross-correlation
of signals}, IEEE Trans. Inform. Theory {\bf 20} (1974) 397-399.




\end{thebibliography}
\end{document}